\documentclass[11pt]{article}
\usepackage{amscd,amssymb,stmaryrd}
\usepackage{amsmath}
\usepackage{graphicx}
\usepackage{subcaption}
\usepackage[utf8]{inputenc}
\usepackage[export]{adjustbox}
\usepackage{wrapfig}
\usepackage{amstext}
\usepackage{pstricks,pst-node,pst-plot,pst-coil}
\usepackage{amsthm}

\usepackage{amsmath}
\usepackage{color}
\usepackage{mathtools}

\usepackage{hyperref}
\usepackage[english]{babel}
\usepackage{booktabs}
\usepackage{mathrsfs}
\usepackage{comment}
\usepackage{float}
\usepackage[linesnumbered,ruled,vlined, ruled,vlined]{algorithm2e}

\theoremstyle{plain}
\theoremstyle{definition}
\newtheorem{theorem}{Theorem}[section]

\usepackage{tikz}
\usepackage{pifont}%

\usepackage{diagbox}
\newcommand{\xmark}{\ding{55}}%
\usetikzlibrary{decorations.pathreplacing,calligraphy}
\usepackage{pgfplots}
\pgfplotsset{compat=1.11}
\usepgfplotslibrary{groupplots}

\title{BelNet: Basis enhanced learning, a mesh-free neural operator}

\date{}

\author{ Zecheng Zhang\footnote{Department of Mathematics, Carnegie Mellon University, Pittsburgh, PA 15213, USA. (Email: zecheng.zhang.math@gmail.com)}, ~Wing Tat Leung \footnote{Department of Mathematics, City University Hong Kong, Hong Kong, China. (Email: wtleung27@cityu.edu.hk)}, 
~Hayden Schaeffer \footnote{Department of Mathematics, UCLA, Los Angeles, CA 90095.} }

\begin{document}

\maketitle

\begin{abstract}
 Operator learning trains a neural network to map functions to functions.  An ideal operator learning framework should be mesh-free in the sense that the training does not require a particular choice of discretization for the input functions, allows for the input and output functions to be on different domains, and is able to have different grids between samples. We propose a mesh-free neural operator for solving parametric partial differential equations.
The basis enhanced learning network (BelNet) projects the input function into a latent space and reconstructs the output functions. In particular, we construct part of the network to learn the ``basis'' functions in the training process. This generalized the networks proposed in \cite{chen1995universal,chen1995approximation} to account for differences in input and output meshes. 
Through several challenging high-contrast and multiscale problems, we show that our approach outperforms other operator learning methods for these tasks and allows for more freedom in the sampling and/or discretization process. 
\end{abstract}

\section{Introduction}
\label{sec_intro}
Data-driven approaches are becoming a competitive and viable means for solving some challenging problems in scientific computing. In particular, machine learning and data science techniques are being developed specifically for scientific tasks such as predicting solutions from complicated spatiotemporal data, approximating trajectories from unknown dynamical systems, or constructing surrogate functions for complex nonlinear operations. One of the underlying principles in these approaches is to either augment or approximate the system using observed data. In particular, these approaches are geared toward solving parametric partial differential equations (PDE). Although one should expect that a rigorous and tailored numerical method for solving a specific PDE will typically outperform a general machine learning (ML) algorithm, one of the key contributions of these new ML approaches is the ability to use less prior information. For example, the ML approaches do not necessarily need to have knowledge of the underlying governing system. However, their ability to correctly capture and model indirect physically relevant terms (e.g. energy or conserved quantities) is still an open problem.

 In this work, we will focus on the task of computing a family of solutions to parametric PDE, i.e. a PDE whose model/equations, boundary, and/or initial data are parametrized.  For example, the conductivity of the heat equation relies on the properties of the underlying material, which can be characterized by some set of parameters. Thus given a set of parameters, one has an equation or algorithm that computes the conductivity and thereby generates a model for the spatiotemporal dynamics. The goal is to find the solution as a function of the parameters, bypassing the need to completely solve the nonlinear system. Recently, deep learning strategies have shown success in this direction \cite{lu2021learning, li2020fourier, zhang2020learning}. 

The standard approaches in ML focus on learning functions, i.e. methods that train a mapping from the parameters to the solution of the equation (at a finite set of locations) based on some pre-generated data. We will refer to these approaches as \textit{function learning}, which includes standard multilayer perceptrons, feedforward neural networks, recurrent neural networks, etc. During the testing phase, one obtains a new set of parameters and the trained network predicts the solution in one forward pass, which can be computationally efficient. However, function learning approaches are restrictive.  First, the neural network may not be easy to learn due to curse-of-dimensionality. Specifically, when training a function that maps the parameters to the solution of a given PDE, the input (i.e. the parameters) and output (i.e. the value of the solution on a discrete set of time and/or space points) will have a very large dimension. Returning to the example, let us consider the 2-dimensional heat equation whose domain is $[0, 1]^2$. If the solution is defined on a spatial mesh with resolution $100\times 100$, then the dimension of the solution vector (i.e. the output) is $10^4$. As a result, the neural network often requires a large set of parameters which can result in a loss of efficiency and training accuracy. Secondly, for function learning, since the output dimension is fixed, all of the output samples must share the same structure. That is, the outputs must be represented on the same mesh or sequence of points. For example, if the majority of observed outputs are defined on a $100\times 100$ mesh, but there are some samples whose solutions are defined on $99\times 99$ mesh or on a different set of $100\times 100$ points, then the standard approaches will fail. This can be very limiting since one does not necessarily have control over the sampling process and may require the system to adapt to different meshes. Lastly, the predicted values must be on the same mesh in order to be consistent with the trained model.

To overcome these issues, \textit{operator learning} \cite{chen1995universal,chen1995approximation} was introduced and has seen success in various scientific computing tasks \cite{lu2021learning, li2020fourier, kovachki2021neural, o2022derivative, lanthaler2022nonlinear, patel2022variationally, ong2022integral, qian2022reduced, de2022cost}. The operator learning approach trains a map between functions rather than fixed-sized vectors, which allows for an accurate representation of the underlying dependencies between the input and output. In particular, operator learning can handle the training of parametric PDEs since the outputs are represented (correctly) by functions rather than fixed dimensional vectors. Formally, let $V$ and $U$ be two function spaces whose domains are $K_1$ and $K_2\subset\mathbb{R}^d$. Let $G$ be the target operator where $G(u)(\cdot)\in U$ for $u\in V$. We want to learn the operator $G$ by the network $G_{\theta}$, i.e.: 
\begin{align}
    G(u)(x) \approx G_{\theta}(u)(x),
    \label{intro_operator}
\end{align}
for $x\in K_2$, where $\theta$ denotes all network parameters. From the structural point of view, the training and expressibility of function learning and operator learning differ. For function learning, the input is the set of parameters that characterize the input function $u$, while in operator learning, the inputs contain both the parameters and the independent variables $x$ used in evaluating the solution. Additionally, in function learning, the output is the solution given a particular discrete representation, but in operator learning, the output is the solution at the point prescribed by the independent variables $x$. Thus with operator learning, one can also predict the solution at other locations within the domain, providing a more flexible representation.

In \cite{chen1995universal}, the authors provided a universal approximation theorem for nonlinear continuous operators. The original application was for modeling solutions to dynamical systems using a particular shallow architecture.  The Deep Operator Neural Network (DON) \cite{lu2021learning}, extended the results of \cite{chen1995universal} (in particular Theorem 5) to account for deep architectures and proposed several new applications. DON has had many extensions and variants. When the governing equations are known, one can incorporate this information using the so-called physics-informed structure \cite{wang2021learning}. In that setting, the user can use fewer samples since the equations provide additional information to the training. Optimization algorithms which exploit DON structure have also been proposed to handle noisy data and train DON \cite{lin2023b, psaros2022uncertainty}, which may be more practical for real data. DON has also been generalized to a wider class of nonlinear approximation problems using shifts \cite{lanthaler2022nonlinear}. Also, detailed numerical experiments on the cost-accuracy trade-off for neural operators are given in \cite{de2022cost}.

The Fourier neural operator (FNO) \cite{li2020fourier, kovachki2021universal} was designed based on a nonlinear generalization of the kernel integral representation for some operators and makes use of the convolutional or Fourier network structure. It was shown to be a universal approximation for operators in \cite{kovachki2021universal}, with respect to the Hilbert space norm (i.e. $H^s$).
FNO has had several extensions with applications to different areas. Some examples include but are not limited to: global weather prediction \cite{pathak2022fourcastnet}, multiphase flow \cite{wen2022u, zhang2022fourier}, and solving PDE with complex geometry \cite{li2022fourier}.
It is worth noting that these neural operator formulations are connected to other approaches in scientific machine learning. FNO can also be seen as a deep generalization to the random Fourier feature approaches \cite{rahimi2007random, rahimi2008weighted} based on approximating kernel integral operators via randomized Fourier networks, see also  \cite{bach2015equivalence,  rudi2017generalization, hashemi2023generalization, chen2021conditioning, chen2022concentration}. In addition, the FNO is connected to other Fourier learning approaches for PDE, including the compressive spectral methods \cite{schaeffer2013sparse, mackey2014compressive}, which trains a sparse Fourier basis representation for the solution of a given PDE; and the sparse optimization framework for PDE discovery \cite{schaeffer2017learning}, which trains a sparse representation of the PDE operator using a mixture of Fourier and polynomial features.

There are several important open problems with operator learning. First, one often needs to discretize the input function $u$ for each sample, i.e. we observe the inputs at some locations or with respect to some basis. However, we do not want the construction of the operator to depend on a fixed choice of locations or basis elements used in the discretization of the inputs since the locations may differ. If a neural operator allows for different input locations in its training, then it is called \textit{discretization-invariant}. Note that the vanilla DON is not discretization-invariant. Although some methods have been proposed to make DON discretization-invariant \cite{li2022learning}, the universal approximation theorem for DON may not support it. Second, we call a neural operator approach \textit{prediction-free} if the discretization mesh for the input function $u$ and the output function $G(u)(\cdot)$ can differ, i.e. we are not restricted to mapping functions on the same input and output grids. Lastly, a method is \textit{domain-independent} if the input function $u$ and its image $G(u)(\cdot)$ can be defined on different domains. For FNO, because of the use of the fast Fourier transform, it is not prediction-free or domain-independent. If a method has all three properties, we call it  \textit{mesh-free}. This work aims to introduce a mesh-free operator learning approach, provide numerical evidence, and show how it can be applied to challenging multiscale problems.

\subsection{Our Contributions}
We propose a mesh-free operator learning approach that builds from the results of \cite{chen1995universal} to incorporate more information on the input and output meshes. Our contributions are as follows. 
\begin{enumerate}
    \item We construct a {mesh-free} operator learning framework called the basis enhanced learning network (BelNet).
    \item We show that BelNet outperforms classical ML approaches and recent neural operators for challenging multiscale systems, i.e. those with high-contrast coefficients or dependencies.
   \item Through several numerical experiments, we show that BelNet achieves consistent results between using a fixed grid or free grid, providing evidence that it generalizes the previous operator learning approaches.
   \end{enumerate}

The paper is organized as follows. In Section \ref{sec_overview}, we review the popular {operator learning} frameworks and introduce the formal definition of the {mesh-free}. We then discuss the details of the basis enhanced learning (BelNet) in Section \ref{sec_bel}. Lastly, Section \ref{sec_numerical} presents some numerical results.

\section{Overview of Operator Networks}
\label{sec_overview}

We review the universal approximation theory and the two recent learning architectures based on approximating operators. We will then introduce the formal definition of {mesh-free}.
 
\subsection{Theory and Architectures}
The classical universal approximation theorems (UAP) for shallow \cite{cybenko1989approximation} and multilayer networks \cite{hornik1989multilayer,pinkus1999approximation} show that for any continuous function, there exists a feedforward neural network with a non-polynomial activation of a certain width (and depth) that can approximate the function with arbitrary accuracy in the sup-norm. The hidden layers' width and depth depend on the accuracy level. Numerical methods for scientific computing and other applications that use the classical UAP are developed for training functions and thus can scale poorly (in terms of trainable parameters or nodes) when the inputs are of very high dimension or have complex structures.

For applications in which the inputs to the model are functions themselves, the classical universal approximation results are insufficient. In \cite{chen1995universal}, it was shown that with a sole hidden layer, one could construct a network that can approximate any operator. These \textit{neural operators} are useful for learning families of solutions to a given PDE since they can approximate mappings that have functional parametric dependencies. We will refer to the universal approximation theorem for operators \cite{chen1995universal} as UAP-O in this work and state it below for completeness.
\begin{theorem}[UAP-O, \cite{chen1995universal}]
    Suppose $g$ is a Tauber-Wiener (TW) function\footnote{If a function $g:\mathbb{R}\rightarrow \mathbb{R}$ (continuous or discontinuous) satisfies that all linear combinations $\Sigma_{i = 1}^Nc_ig(\lambda_i x+\theta_i)$ are dense in $C[a, b]$, where $c_i, \lambda_i, \theta_i\in\mathbb{R}$, then $g$ is called a Tauber-Wiener (TW) function.}, $X$ is Banach, $K_1\subset X$ and $K_2\subset \mathbb{R}^d$ are compact. Let $V\subset C(K_1)$ be compact and $G: V\rightarrow C(K_2)$ be a nonlinear continuous operator.
    Then for any $\epsilon>0$, there are positive integers $M, N, K$, constants $c_i^k, \zeta_k, \theta_i^k, \varepsilon_{ij}^k\in\mathbb{R}$, points $\omega_k\in\mathbb{R}^d$, $y_j\in K_1$, $i = 1, ..., M$, $k = 1, ..., K$, $j = 1, ...., N$ such that
    \begin{align*}
        \bigg| G(u)(x) - \sum_{k = 1}^K \sum_{i = 1}^M c_i^k\, g\left(\sum_{j = 1}^N\varepsilon_{ij}^ku(y_j)+\theta_i^k\right)\, g(\omega_k\cdot x+\zeta_k)\bigg|<\epsilon
    \end{align*}
    holds for all $u\in V$ and $x\in K_2$.
\end{theorem}
The DON framework \cite{lu2021learning, lu2022comprehensive, lanthaler2022error} extends the UAP-O to continuous functions and then builds the network based on the following general UAP theorem stated below.
\begin{theorem}[DON Approximation \cite{lu2021learning, lu2022comprehensive, lanthaler2022error}]
    Suppose $X$ is Banach, $K_1\subset X$ and $K_2\subset \mathbb{R}^d$ are compact. Let $V\subset C(K_1)$ be compact and $G: V\rightarrow C(K_2)$ be a nonlinear continuous operator. Then for any $\epsilon>0$, there exists positive integers $N, K$, continuous function $f:\mathbb{R}^d\rightarrow\mathbb{R}^K$, $g:\mathbb{R}^N\rightarrow\mathbb{R}^K$ and $y_1, ..., y_N\in K_1 $ such that
    \begin{align*}
        \left| G(u)(x) - \langle\underbrace{g\big(u(y_1), ..., u(y_N)\big)}_{branch}, \underbrace{f\big(x\big)}_{trunk}\rangle\right|<\epsilon,
    \end{align*}
    for any $u\in V$ and $x\in K_2$. Here $\langle\cdot, \cdot\rangle$ denotes the inner product in $\mathbb{R}^K$.
    \label{deeponet_uap}
\end{theorem}
Theorem~\ref{deeponet_uap} has been used as the basis for many recent developments in neural network approximations based on the DON framework \cite{lu2022comprehensive, lanthaler2022error}. Note that, from the structure of the network used in Theorem~\ref{deeponet_uap}, we see that the collection of training sensors $\{y_j\}_{j=1}^N$ are ``fixed'' for all the sample $u\in V$, i.e. the locations of the observations are the same for each sample. For example, to train the solution operator for nonlinear PDE, the grid (i.e. the choice of discretization) must be the same for all input functions in order of  Theorem~\ref{deeponet_uap} to hold. 

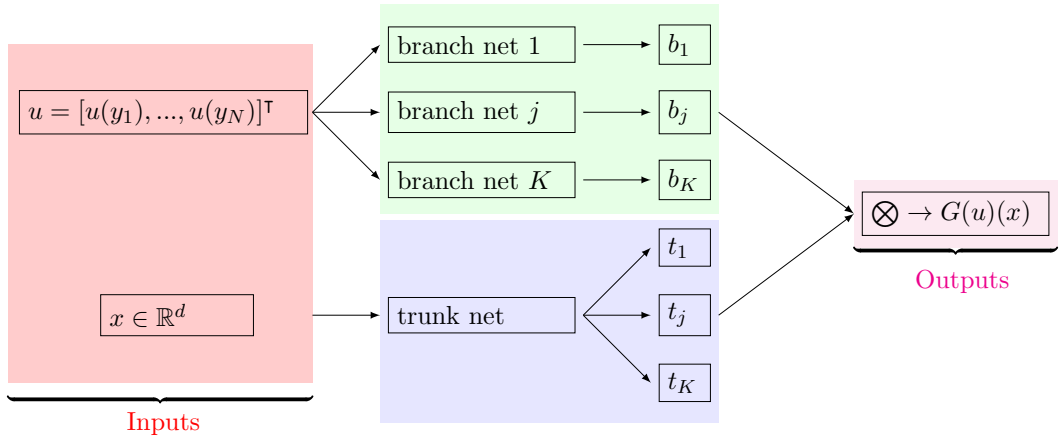
\begin{figure}[H]
\centering
\scalebox{.9}{\begin{tikzpicture}[scale = 1]
 \fill [green!10] (2, 4.5) rectangle (7, 7.6);
  \fill [blue!10] (2, 1.4) rectangle (7, 4.4);
    \fill [red!20] (-3.5, 2) rectangle (1, 7);

    \draw[ultra thick] [decorate,
    decoration = {calligraphic brace, mirror}] (-3.5, 1.8) --  (1, 1.8);
\node at (-1.2, 1.4) {\textcolor{red}{Inputs}};

\fill [magenta!10] (9, 4) rectangle (12, 5);
    \draw[ultra thick] [decorate,
    decoration = {calligraphic brace, mirror}] (9, 4) --  (12, 4);
\node at (10.6, 3.5) {\textcolor{magenta}{Outputs}};

\node[draw, text width=4cm] at (-1.2, 6) {$u = [u(y_1), ..., u(y_N)]^\intercal$};

 \draw [-latex ](1,6) -- (2, 7);
 \node[draw, text width = 2.5cm] at (3.5, 7) {branch net $1$};
 \draw [-latex ](5, 7) -- (6, 7);
 \node[draw, text width = 0.5cm] at (6.5, 7) {$b_1$};
 
 \draw [-latex ](1,6) -- (2, 5);
  \node[draw, text width = 2.5cm] at (3.5, 5) {branch net $K$};
  \draw [-latex ](5, 5) -- (6, 5);
  \node[draw, text width = 0.5cm] at (6.5, 5) {$b_K$};
  
  \draw [-latex ](1,6) -- (2, 6);
 \node[draw, text width=2.5cm] at (3.5, 6) {branch net $j$};
 \draw [-latex ](5, 6) -- (6, 6);
\node[draw, text width = 0.5cm] at (6.5, 6) {$b_j$};

 \node[draw, text width = 2cm] at (-1, 3) {$x\in\mathbb{R}^d$};
  \draw [-latex ](1, 3) -- (2, 3);
 \node[draw, text width = 2.5cm] at (3.5, 3) {trunk net};
 
 \draw [-latex ](5, 3) -- (6, 4);
\node[draw, text width = 0.5cm] at (6.5, 4) {$t_1$};

  \draw [-latex ](5, 3) -- (6, 3);
\node[draw, text width = 0.5cm] at (6.5, 3) {$t_j$};

   \draw [-latex ](5, 3) -- (6, 2);
 \node[draw, text width = 0.5cm] at (6.5, 2) {$t_K$};

\draw [-latex ](7, 6) -- (9, 4.5);
  \draw [-latex ](7, 3) -- (9, 4.5);

  \node[draw, text width = 2.5cm] at (10.5, 4.5) {$\bigotimes\rightarrow G(u)(x)$};

\end{tikzpicture}}
\caption{Stacked version DON. $\bigotimes$ denotes the inner product in $\mathbb{R}^K$.}
\label{deepo_structure}
\end{figure}

Figure \ref{deepo_structure} presents an illustration of the architecture of DON based on Theorem \ref{deeponet_uap}. In particular, note that DON depends on the independent variable $x\in\mathbb{R}^d$ within the separate trunk subnetwork (highlighted in blue in Figure \ref{deepo_structure}). As a result, the input and output domains can differ, and the choice of discretization for the two domains do not need to be the same. However, the set of sensor locations $\{y_j\}_{j=1}^N$ is fixed.

Another parallel neural operator introduced in \cite{li2020fourier} is the Fourier neural operator (FNO), which was originally designed as an approximation to the kernel integral operator \cite{li2020fourier, kovachki2021universal} for training a family of solutions to PDE. FNO takes the following form \cite{kovachki2021neural}:
\begin{align}
    \mathcal{FNO} (u) = Q\circ \mathcal{L}_L\circ...\circ\mathcal{L}_1\circ R(u)
    \label{intro_fno}
\end{align}
where $Q$ and $R$ are the linear projection and lifting operators, $\mathcal{L}_l$ is the neural (sub)network of the form:
\begin{align*}
    \mathcal{L}_l(u)(x) = \sigma\bigg(W_lu(x)+b_l(x)+\mathcal{F}^{-1}\big(P_l(k)\cdot \mathcal{F}(u)(k)\big)(x) \bigg),
\end{align*}
where $W_l$ and $b_l$ are the weight and bias.
In this formulation, $\mathcal{F}$ and $\mathcal{F}^{-1}$ are the Fourier transformation and its inverse, and thus the majority of hidden parameters are trained in the Fourier domain. In particular, $P_l$ defines the coefficients of a non-local, linear mapping in the Fourier domain. Although graph-based formulations are also possible \cite{li2020neural}, writing the network as the composition of nonlinear operations in the Fourier domain provides several computational benefits, including lowering the computational cost and simplifying the layer structure. In addition, \textcolor{black}{FNO is discretization invariant in the sense that the network is invariant to the different meshes (see Theorem 5 of \cite{kovachki2021universal}). However, the input and output functions must be defined in the same domain on the same grid.}

Note that in both neural operator frameworks, the sensor locations are fixed. In applications where the sensors are moving (like in particle tracking) or in experiments where not all sensors provide useful information, we would like to be able to use different samples $\{y_j\}_{j=1}^N$. In particular, we construct a neural operator based on DON that is mesh-free, i.e. independent of the locations used in the discretization, prediction, or domains. 

\subsection{Discretization-Invariant}
\begin{figure}[H]
\centering
\begin{tikzpicture}[scale=1]

\draw[step=0.5cm, gray] (-4,0) grid (-2, 2);

\foreach \i in {0, ..., 4} {
        \foreach \j in {0, ..., 4}{
\filldraw[red] (-4+\i/2, \j/2) circle (2pt) node[anchor=west]{};
        }
        }
\node at (-3.5, 2.3) {\tiny Sample 1};   

\draw[step = 0.5cm, gray] (-1.5-0.001, 0) grid (0.5, 2);

\foreach \i in {0, ..., 4} {
        \foreach \j in {0, ..., 4}{
\filldraw[red] (-1.5+\i/2, \j/2) circle (2pt) node[anchor=west]{};
        }
        }
\node at (-1, 2.3) {\tiny Sample 2}; 

\draw[ultra thick] [decorate,
    decoration = {calligraphic brace, mirror}] (-4, -0.5) --  (0.5, -0.5);
\node[text width = 4.5cm] at (-1.5, -1.8) {An uniform mesh to evaluate the input functions (eg: finite difference methods). };

\draw [draw=black] (1, 0) rectangle (3, 2);
\filldraw[red] (1.5, 1.5) circle (2pt) node[anchor=west]{};
\filldraw[red] (2, 1) circle (2pt) node[anchor=west]{};
\filldraw[red] (2.5, 0.5) circle (2pt) node[anchor=west]{};

\node at (1.5, 2.3) {\tiny Sample 1};  
 \draw [draw=black] (3.5, 0) rectangle (5.5, 2);
\filldraw[red] (1.5+2.5, 1.5) circle (2pt) node[anchor=west]{};
\filldraw[red] (2+2.5, 1) circle (2pt) node[anchor=west]{};
\filldraw[red] (2.5+2.5, 0.5) circle (2pt) node[anchor=west]{};
\node at (4, 2.3) {\tiny Sample 2}; 

\draw[ultra thick] [decorate,
    decoration = {calligraphic brace, mirror}] (1, -0.5) --  (5.5, -0.5);
\node[text width = 4.5cm] at (3.5, -2) {{Discretization-dependent}: input functions of two samples must be sampled at the same sensors (red dots)};

\draw [draw=black] (6, 0) rectangle (8, 2);
\filldraw[red] (1.5+2.5*2, 1.5) circle (2pt) node[anchor=west]{};
\filldraw[red] (2+2.5*2, 1) circle (2pt) node[anchor=west]{};
\filldraw[red] (2.5+2.5*2, 0.5) circle (2pt) node[anchor=west]{};
\node at (6.5, 2.3) {\tiny Sample 1};  

\draw [draw=black] (8.5, 0) rectangle (10.5, 2);
\filldraw[red] (1.5+2.5*3, 0.5) circle (2pt) node[anchor=west]{};
\filldraw[red] (2+2.5*3, 1.5) circle (2pt) node[anchor=west]{};
\filldraw[red] (2.5+2.5*3, 1) circle (2pt) node[anchor=west]{};
\node at (9, 2.3) {\tiny Sample 2}; 

\draw[ultra thick] [decorate,
    decoration = {calligraphic brace, mirror}] (6, -0.5) --  (10.5, -0.5);
\node[text width = 4.5cm] at (8.5, -2) {{Discretization-invariant}: input functions of two samples can be sampled at the \textit{different} sensors (red dots)};

\end{tikzpicture}
\caption{Illustration of discretization-invariant. Discretization-invariant (right two images) does not require different samples to share the identical discretization of the parameter space, i.e. parameters can be evaluated at different sets of points (dots). }
\label{intro_dis_free}
\end{figure}
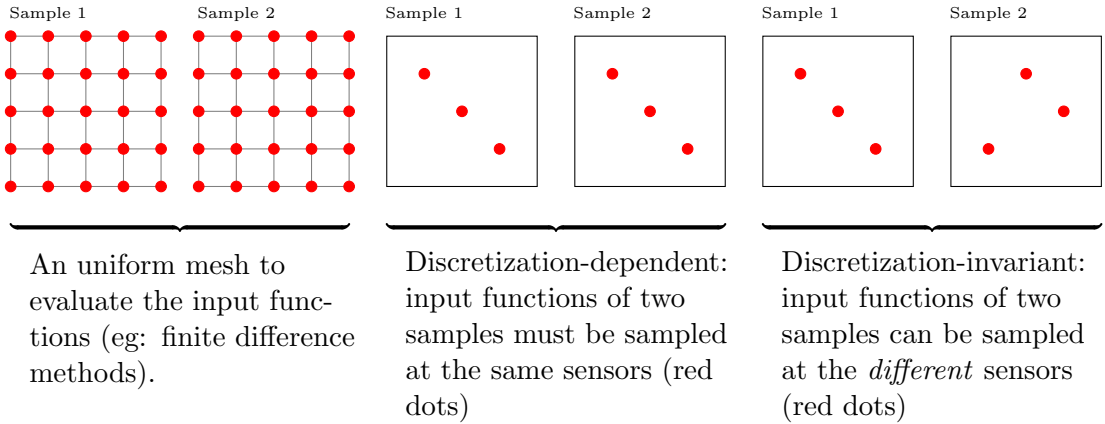

The goal is to construct an operator $G(u)(x): V\rightarrow U$, for $u\in V$, $G(u)(x)\in U$, where $V$ and $U$ are two function spaces.
The function $u$ (which can be the controls, boundary conditions, initial states, etc.) is an input to the network and will only be observed on a set of sensor locations. This is represented by a particular choice of discretization. In practice, the discretization of $u$ could be determined by the data or the sensor locations and thus may not be fixed between different experiments. However, based on the UAP-O, $u$ must be observed at some sensor points to obtain the approximation guarantees. We call a neural operator {discretization-invariant} if the approximation does not require the locations of the sensors $\{y_j\}_{j=1}^N$ to be fixed. As an illustrative example, Figure \ref{intro_dis_free} considers three cases between two different samples, where the sample refers to the different input functions $u$. The first plot of Figure \ref{intro_dis_free} shows the uniform grid used in the construction of finite differences or spectral methods, which require each sampled input, e.g. sample 1 $\{u_1(y_j)\}_{j=1}^N$,  and sample 2, $\{u_2(y_j)\}_{j=1}^N$, share the same grid. This is also the case in the middle plot, which uses the same three sensor locations. Our approach is discretization-invariant in the sense that the sensors can be distributed independently for different samples displayed in the plot on the right in Figure \ref{intro_dis_free}. Based on the construction, the FNO approach is discretization-invariant since \textcolor{black}{the layers are trained in the Fourier domain and are thus evaluated in the spatial domain via a projection using the Fourier basis elements}. The vanilla DON formulation is not discretization-invariant. Some recent extensions of DON \cite{li2022learning} could potentially make it discretization-invariant but showing this holds remains to be proven.

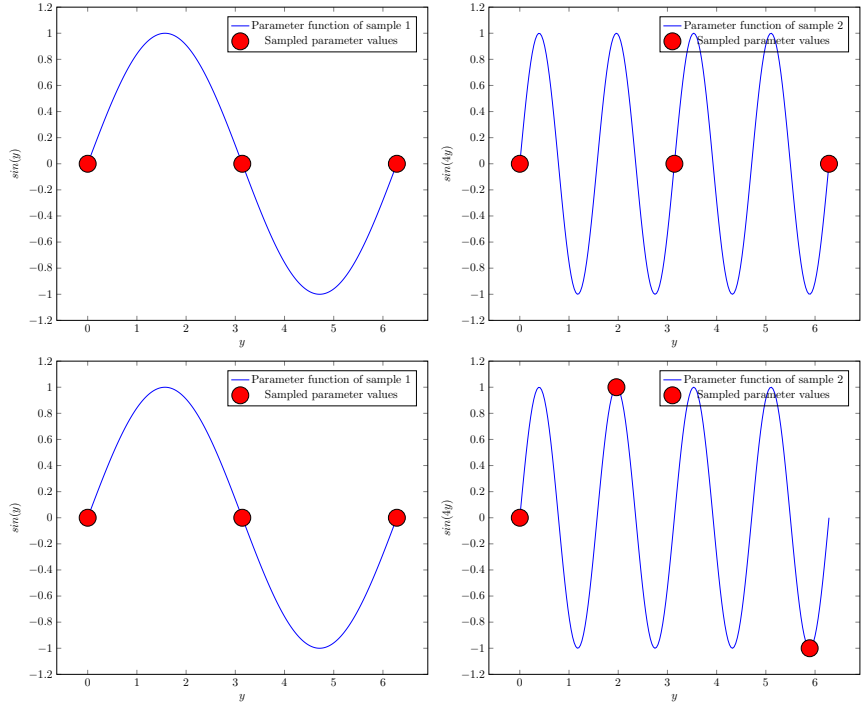
\begin{figure}[h]
\centering
\begin{tikzpicture}[scale = 0.4]
\begin{axis}[domain=0:2*pi, legend pos= north east, ylabel = $sin(y)$, xlabel = $y$, height = 12cm]
  \addplot[blue, samples = 200] {sin( deg(x)  ) }; 
  \addplot[only marks,  mark=*, mark size=8pt, mark options={fill=red}, samples = 3] {sin( deg(x) )}; 
  \legend{Parameter function of sample 1, Sampled parameter values}
  \node[text width = 6cm] at (3, 6) {Sample 1 parameter function};
\end{axis}
\end{tikzpicture}
\begin{tikzpicture}[scale = 0.4]
\begin{axis}[domain=0:2*pi, legend pos= north east, legend style={fill=none}, ylabel = $sin(4y)$, xlabel = $y$, , height = 12cm]
  \addplot[blue, samples = 200] {sin( deg(4*x)  ) }; 
  \addplot[only marks,  mark=*, mark size=8pt, mark options={fill=red}, samples = 3] {sin( deg(4*x) )}; 
  \legend{Parameter function of sample 2, Sampled parameter values};
\end{axis}
\end{tikzpicture}
\begin{tikzpicture}[scale = 0.4]
\begin{axis}[domain=0:2*pi, legend pos= north east, ylabel = $sin(y)$, xlabel = $y$, , height = 12cm]
  \addplot[blue, samples = 200] {sin( deg(x)  ) }; 
  \addplot[only marks,  mark=*, mark size=8pt, mark options={fill=red}, samples = 3] {sin( deg(x) )}; 
  \legend{ Parameter function of sample 1,  Sampled parameter values}
  \node[text width = 6cm] at (3, 6) {Sample 1 parameter function};
\end{axis}
\end{tikzpicture}
\begin{tikzpicture}[scale = 0.4]
\begin{axis}[domain=0:2*pi, legend pos= north east, legend style={fill=none}, ylabel = $sin(4y)$, xlabel = $y$, height = 12cm]
  \addplot[blue, samples = 200] {sin( deg(4*x)  ) }; 
  \addplot [only marks, mark=*, mark size=8pt, mark options = {fill = red}] coordinates { (0,0) (2*pi-pi/8, -1 ) (pi-pi/8-pi/4, 1 ) };
  \legend{ Parameter function of sample 2, Sampled parameter values};
\end{axis}
\end{tikzpicture}
\caption{Illustration of discretization-invariant. Discretization-invariant (lower two images) can allow sampling at different sensors in the domain for two input functions, such that the sampled function values are different. However, not being discretization-invariant (top two images) may result in the sampled function values being the same for two different samples. }
\label{intro_dis_free2}
\end{figure}

A benefit to discretization-invariance is the potential for avoiding aliasing between samples, which is a particular issue for the sampled parametric dependencies in multi-scale problems. This can impact both the network training and the overall accuracy. We illustrate this idea in Figure \ref{intro_dis_free2}.

\subsection{{Predication-Free} and {Domain-Independent}}
\begin{figure}[H]
\centering
\begin{tikzpicture}[scale=1]

\draw[step=0.5cm, gray] (-4,0) grid (-2, 2);

\foreach \i in {0, ..., 4} {
        \foreach \j in {0, ..., 4}{
\filldraw[red] (-4+\i/2, \j/2) circle (2pt) node[anchor=west]{};
        }
        }
\node at (-3., 2.3) {\tiny Input mesh};   

\draw[step = 0.5cm, gray] (-1.5-0.001, 0) grid (0.5, 2);

\foreach \i in {0, ..., 4} {
        \foreach \j in {0, ..., 4}{
\filldraw[blue] (-1.5+\i/2, \j/2) circle (2pt) node[anchor=west]{};
        }
        }
\node at (-1+0.5, 2.3) {\tiny Output mesh}; 

\draw[ultra thick] [decorate,
    decoration = {calligraphic brace, mirror}] (-4, -0.5) --  (0.5, -0.5);
\node[text width = 4.5cm] at (-1.5, -1.55) {An uniform mesh to sample the input function (eg: finite difference methods). };

\draw [draw=black] (1, 0) rectangle (3, 2);
\filldraw[red] (1.5, 1.5) circle (2pt) node[anchor=west]{};
\filldraw[red] (2, 1) circle (2pt) node[anchor=west]{};
\filldraw[red] (2.5, 0.5) circle (2pt) node[anchor=west]{};

\node at (1.5+0.5, 2.3) {\tiny Input mesh};  
 \draw [draw=black] (3.5, 0) rectangle (5.5, 2);
\filldraw[blue] (1.5+2.5, 1.5) circle (2pt) node[anchor=west]{};
\filldraw[blue] (2+2.5, 1) circle (2pt) node[anchor=west]{};
\filldraw[blue] (2.5+2.5, 0.5) circle (2pt) node[anchor=west]{};
\node at (4+0.5, 2.3) {\tiny Output mesh}; 

\draw[ultra thick] [decorate,
    decoration = {calligraphic brace, mirror}] (1, -0.5) --  (5.5, -0.5);
\node[text width = 4.5cm] at (3.5, -2) {Not {prediction-free}: output functions must be evaluated at the same points as the input functions.};

\draw [draw=black] (6, 0) rectangle (8, 2);
\filldraw[red] (1.5+2.5*2, 1.5) circle (2pt) node[anchor=west]{};
\filldraw[red] (2+2.5*2, 1) circle (2pt) node[anchor=west]{};
\filldraw[red] (2.5+2.5*2, 0.5) circle (2pt) node[anchor=west]{};
\node at (6.5+0.5, 2.3) {\tiny Input mesh};  

\draw [draw=black] (8.5, 0) rectangle (10.5, 2);
\filldraw[blue] (1.5+2.5*3, 0.5) circle (2pt) node[anchor=west]{};
\filldraw[blue] (2+2.5*3, 1.5) circle (2pt) node[anchor=west]{};
\filldraw[blue] (2.5+2.5*3, 1) circle (2pt) node[anchor=west]{};
\node at (9+0.5, 2.3) {\tiny Output mesh}; 

\draw[ultra thick] [decorate,
    decoration = {calligraphic brace, mirror}] (6, -0.5) --  (10.5, -0.5);
\node[text width = 4.5cm] at (8.5, -2) {{Prediction-free}: output can be evaluated at any point which is independent of the input discretization.};

\end{tikzpicture}
\caption{Illustration of {prediction-free}. }
\label{intro_pred_free}
\end{figure}
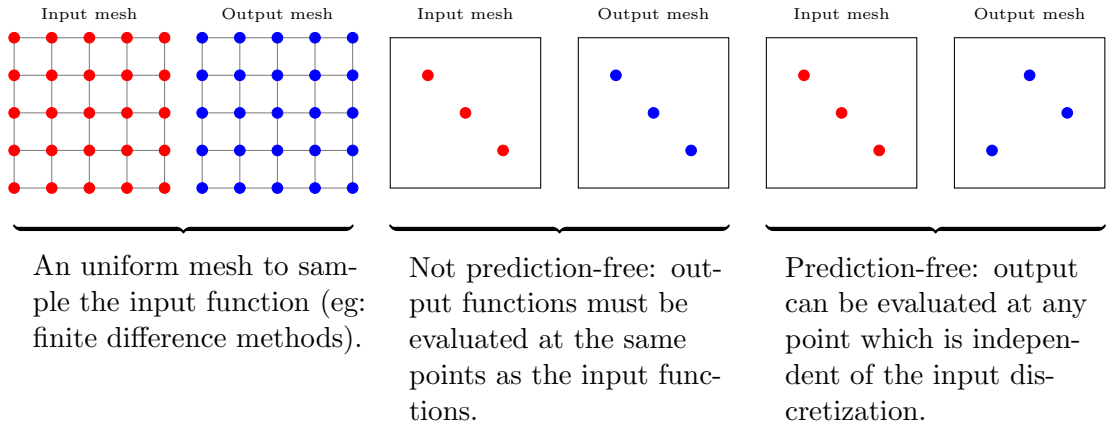

The goal is to build the neural operator $G(u)$ such that one can evaluate the function values $G(u)(x)$ for $x\in K_2$, i.e. for interpolation or prediction. We say that the neural operator is {prediction-free} if the  evaluation point $x$ can be any point in $K_2$ (output function is free of discretization), independent of the choice of the discretization of the input function used in the training of the neural operator. Figure \ref{intro_pred_free} provides an example of {prediction-free} (the mesh on the right) as opposed to the first two plots, which are not {prediction-free} in the sense that the output mesh depends on the choice of the input mesh. 
By construction, the DON is a {prediction-free} approximation based on UAP-O. Vanilla FNO is not {prediction-free} since the output mesh must match the discretization of the input mesh \cite{lu2022comprehensive}. While requiring a neural operator to be prediction-free may be challenging, it can be shown to have several benefits. One advantage is that it leverages the operator learning principle and overcomes the limitations of function learning, e.g. lower network complexity, less dependence on data sampling process, and testing flexibility (see also Section \ref{sec_intro}).

\begin{figure}[H]
\centering
\begin{tikzpicture}[scale=1.2]

\draw [draw=black] (1, 0) rectangle (3, 2);
\filldraw[red] (1.5, 1.5) circle (2pt) node[anchor=west]{};
\filldraw[red] (2, 1) circle (2pt) node[anchor=west]{};
\filldraw[red] (2.5, 0.5) circle (2pt) node[anchor=west]{};

\node at (1.5+0.5, 2.3) {\tiny Input domain};  
 \draw [draw=black] (3.5, 0) rectangle (5.5, 2);
\filldraw[blue] (1.5+2.5, 1.5) circle (2pt) node[anchor=west]{};
\filldraw[blue] (2+2.5, 1) circle (2pt) node[anchor=west]{};
\filldraw[blue] (2.5+2.5, 0.5) circle (2pt) node[anchor=west]{};
\node at (4+0.5, 2.3) {\tiny Output domain}; 

\draw[ultra thick] [decorate,
    decoration = {calligraphic brace, mirror}] (1, -0.5) --  (5.5, -0.5);
\node[text width = 4.5cm] at (3.5, -1.6) {{Domain-dependent}: input and output functions must have the same domain.};

\draw (6, 1) .. controls (7,2) and (6.5, 1) .. (8,1);
\filldraw[red] (6, 1) circle (2pt) node[anchor=west]{};
\filldraw[red] (6.5, 1.40) circle (2pt) node[anchor=west]{};
\filldraw[red] (8, 1) circle (2pt) node[anchor=west]{};
\node at (6.5+0.5, 2.3) {\tiny Input domain};  

\draw [draw=black] (8.5, 0) rectangle (10.5, 2);
\filldraw[blue] (1.5+2.5*3, 0.5) circle (2pt) node[anchor=west]{};
\filldraw[blue] (2+2.5*3, 1.5) circle (2pt) node[anchor=west]{};
\filldraw[blue] (2.5+2.5*3, 1) circle (2pt) node[anchor=west]{};
\node at (9+0.5, 2.3) {\tiny Output domain}; 

\draw[ultra thick] [decorate,
    decoration = {calligraphic brace, mirror}] (6, -0.5) --  (10.5, -0.5);
\node[text width = 4.5cm] at (8.5, -1.76) {{Domain-independent}: input and output functions can be defined on the different domains. };

\end{tikzpicture}
\caption{Illustration of domain-independence. }
\label{intro_domain_free}
\end{figure}
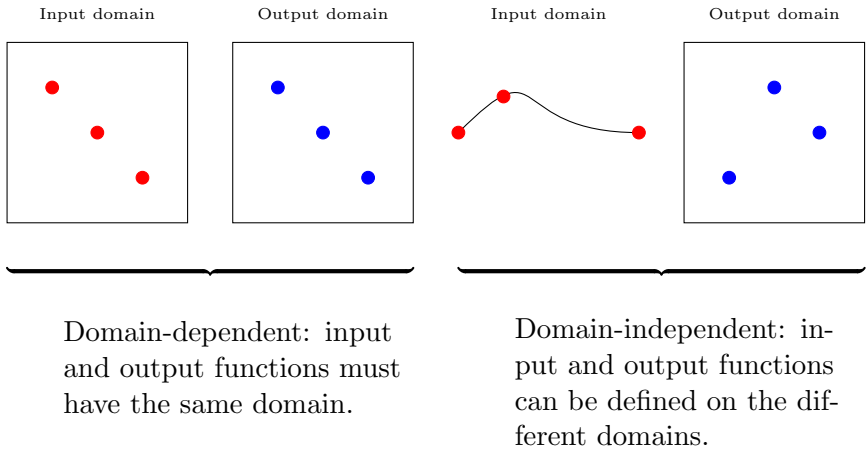

In addition, FNO requires that the input and output function spaces be defined on the same domain, so it is {domain-dependent}. We call an operator learning framework {domain-independent} if the output function domain is independent of the input function domain. As an example of domain-independent, see Figure \ref{intro_domain_free}. Domain-independent methods provide the flexibility to solve problems from broader areas. One example is the parametric time-dependent PDE. While the parametrized input function may only have the spatial dependence, the solution operator relies on both space and time. 

\vspace{0.5em}
\subsection{Mesh-Free}\label{goals}
A flexible neural operator tool is one that is {mesh-free}, that is, a technique that is {discretization-invariant}, {domain-independent} and {prediction-free}. As a comparison, we summarize the comparison between the neural operators in Table \ref{intro_table}, see also \cite{lu2022comprehensive}. We will show that the Basis Enhanced Learning Network (BelNet) is a {mesh-free} extension of the current state-of-the-art neural operators.

\begin{table}[h]
\centering
\begin{tabular}{| c | c |  c | c|}
\hline
& Discretization-invariant & Prediction-free & Domain-independent\\
\hline
DON & \textcolor{red}{\xmark} & \textcolor{green}{\checkmark} & \textcolor{green}{\checkmark} \\
\hline
FNO & \textcolor{green}{\checkmark} & \textcolor{red}{\xmark} & \textcolor{red}{\xmark}\\
\hline
BelNet & \textcolor{green}{\checkmark} & \textcolor{green}{\checkmark} & \textcolor{green}{\checkmark} \\
\hline
\end{tabular}
\caption{Comparison of DON, FNO and our BelNet.}
\label{intro_table}
\end{table}

\section{BelNet Formulation}
\label{sec_bel}
Similar to FNO, BelNet is derived as a nonlinear generalization of a linear kernel operator. 

\textbf{Linear Setting}: We want to approximate an operator $G: V\rightarrow U$ as discussed in Equation (\ref{intro_operator}). As a motivating example, assume that the underlying kernel $\kappa(x, y)$ can be expanded in the following way: $\kappa(x, y) = \sum_{k = 1}^Kp_k(y) q_k(x)$, then
\begin{align*}
    G(u)(x) = \int \kappa(x, y)u(y)dy  = \sum_{k = 1}^K q_k(x)\int p_k(y)u(y)dy.
\end{align*}
Note that this expansion of the kernel can hold in a variety of useful settings. For example, if the kernel is positive definite then there exists a feature map $\phi$ such that $\kappa(x, y)=\langle \phi(y), \phi(x) \rangle$ by Mercer's Theorem. If the kernel is low-rank (or, more generally, if we approximate it by a low-rank expansion) then we can write the kernel in the form $\kappa(x, y) = \sum_{k = 1}^Kp_k(y) q_k(x)$ where $p_k$ and $q_k$ are the right and left singular functions rescaled by the singular values. To approximate the integral, we introduce a quadrature rule:
\begin{align}
    G(u)(x) = \sum_{k = 1}^K q_k(x)\int p_k(y)u(y)dy \approx \sum_{k = 1}^K q_k(x)\sum_{j = 1}^N p_{kj}u_{j},  \label{intro_formulation2}
\end{align}
where $u_j = u(y_j)$, for $y:= [y_i, ..., y_N]^\intercal\subset K_1$, and $p_{kj}\in\mathbb{R}$ for all $j, k$ are the weights defined by $p_k(y_j)$ rescaled by the quadrature coefficients. Thus the natural neural operator in the linear setting becomes
\begin{align}
    G_\theta(u)(x) = \sum_{k = 1}^K q_k(x) \, \sum_{j = 1}^N p_{kj}u_{j}= \sum_{k = 1}^K q_k(x) \, p_{k}^\intercal u,
      \label{form_linear}
\end{align}
where we can train each function $q_k(x)$ by a standard MLP and collect all of the weights used in the MLP's along with $p_{kj}$'s into a vector of trainable parameters $\theta$.

\textbf{Nonlinear Setting}: To generalize Equation~\eqref{form_linear} to a nonlinear operator, we introduce a nonlinear function $g$:
\begin{align*}
    G(u)(x) \approx \sum_{k = 1}^K q_k(x)\, g\left(\sum_{j = 1}^N p_{kj}u_{j}\right), 
     \end{align*}
which incorporates nonlinear dependencies on the points $x$ through the basis functions $q_k$ and nonlinear dependencies in the function value $u_j$. We construct a neural operator $G_{\theta}$ to approximate by training the functions $q_k$ and $g$ and the values $p_{kj}$.  
In particular, we introduce weights and biases, $q^{k}\in\mathbb{R}^{d}$, $W_y^{1, k}\in\mathbb{R}^{N_1\times N}$, $W_y^{2, k}\in\mathbb{R}^{N\times N_1}$, 
$b_x^k\in \mathbb{R}$, $b_y^k\in\mathbb{R}^{N_1}$, where $k = 1, ..., K$, and activation functions $a_x$, $a_y$ and $a_u$ 
such that,
\begin{align}
    G(u)(x) \approx G_{\theta}(u)(x) = \sum_{k = 1}^K  a_x  \left((q^{k})^\intercal x+ b_x^k \right) \, a_u \left(\hat{u}^\intercal W_{y}^{2, k}\left(a_y(W_{y}^{1, k} y + b_y^k )  \right) \right),
    \label{intro_formulation3}
\end{align}
for $x\in K_2\subset \mathbb{R}^d$, $u\in V$,  and where $\hat{u} = [u(y_1), ..., u(y_N)]^\intercal\in\mathbb{R}^N$.

The motivation for this formulation is its {discretization-invariant} nature. To see this, let us assume $V \subset \mathbb{P}_N$, $p_k(y)>0$ and $p_k(y)\in L^1(\Omega)$ where $\Omega$ is the domain. 
For any set of distinct points $\{y_j\}_{j = 0}^N\subset I$, we can define the associated Lagrange basis as,
\begin{align*}
    h_j(y) = \prod_{l=0, k\neq j}^N\frac{y-y_j}{y_j-y_l}, 
\end{align*}
for $j\in[0, N]$. 
It is known that \cite{shen2011spectral},
\begin{align*}
    \int_I p_k(y)u(y)dy = \sum_{j = 0}^N u(y_j)w_{kj}
\end{align*}
where,
$w_{kj} = \int_I h_j(y)p_k(y)dy, j\in[0, N].$
That is, the approximation of the integral in Equation \ref{intro_formulation2} is exact. However, 
it should be noted that $\{y_j\}$ are the sensors of the input function and arbitrary, indicating that the method is independent of the discretization.

\begin{figure}[H]
\centering
\scalebox{.85}{\begin{tikzpicture}[scale = 1]

\fill [red!20] (-3, 0) rectangle (0.9, 7);
\draw[ultra thick] [decorate,
    decoration = {calligraphic brace, mirror}] (-2.5, -0.1) --  (0.5, -0.1);
\node at (-1, -0.5) {\textcolor{red}{Inputs}};

\fill [green!10] (1.8, 4.3) rectangle (7.5, 7.6);
\draw[ultra thick] [decorate,
    decoration = {calligraphic brace, mirror}] (3, 4.2) --  (6, 4.2);
\node at (4.5, 3.9) {\textcolor{green}{Projection nets}};

\fill [blue!10] (3, 0) rectangle (6, 2);

\draw[ultra thick] [decorate,
    decoration = {calligraphic brace, mirror}] (3, -0.1) --  (6, -0.1);
\node at (4.5, -0.5) {\textcolor{blue}{Construction net}};

\fill [magenta!10] (10., 0) rectangle (14.5, 2);

\draw[ultra thick] [decorate,
    decoration = {calligraphic brace, mirror}] (10.8, -0.1) --  (13.8, -0.1);
\node at (12.2, -0.5) {\textcolor{magenta}{Outputs}};
    
 \node[draw, text width = 3.5cm] at (-1, 6) {$y = [y_1, y_2, ..., y_N]^\intercal$};
 
 \draw [-latex ](1,6) -- (2, 7);
 \node[draw, text width = 5cm] at (4.7, 7) {$W_y^{2, 1}a_y(W_y^{1,1}y+b_y^1):=p^1$};
 
 \draw [-latex ](1,6) -- (2, 5);
  \node[draw, text width = 5cm] at (4.7, 5) {$W_y^{2, K}a_y(W_y^{1, K}y+b_y^K):= p^K$};
  
  \draw [-latex ](1,6) -- (2, 6);
 \node[text width=5cm] at (5, 6) {$... ...$ \quad $... ...$ \quad $... ...$};

 \draw [-latex ](7.5, 7) -- (8.5, 6);
 \draw [-latex ](7.5, 6) -- (8.5, 6);
 \draw [-latex ](7.5, 5) -- (8.5, 6);

 \node[draw, text width = 5cm] at (11.5, 6) {Concatenate the outputs $\{p^1, ..., p^K\}$ to get a matrix $P \in\mathbb{R}^{K\times N}$};

\node[draw, text width=3.5cm] at (-1, 3) {$u = [u_1, u_2, ..., u_N]^\intercal$};
\draw [-latex ](1, 3) -- (10, 3);

\node[draw, text width = 1.8 cm, align = center] at (11.5, 3) {$a_u(Pu)$};

\draw [-latex ](11.5, 5) -- (11.5, 3.5);

 \node[draw, text width = 3.5cm] at (-1, 1) {$x = [x_1, x_2, ..., x_d]^\intercal$};
 \draw [-latex ](1,1) -- (3, 1);
\node[draw] at (4.5, 1) {$a_x(Qx + b_x)$};
 \draw [-latex ](6,1) -- (10, 1);

\node[draw, text width = 3.5 cm, align = center] at (12.3, 1) {$[a_u(Pu)]^\intercal a_x(Qx + b_x)$};

\draw [-latex ](11.5, 2.5) -- (11.5, 2.0);

\end{tikzpicture}}
\caption{Illustration of the network structure. Projection nets are $K$ independent fully connected neural network with weights and bias $W_y^{2, k}\in\mathbb{R}^{N\times N_1}$, $W_y^{1, k}\in\mathbb{R}^{N_1\times N}$ and $b_y^k\in\mathbb{R}^{N_1}$. Construction net is a fully connected neural network with weights and bias $Q\in\mathbb{R}^{K\times d}$ and $b_x\in\mathbb{R}^d$. Here $Q = [q^1, q^2, ..., q^K]$, where $q^i$ are defined in Equation (\ref{intro_formulation3}). In addition, $a_x, a_y, a_u$ are activation functions. }
\label{fig_network_s1}
\end{figure}
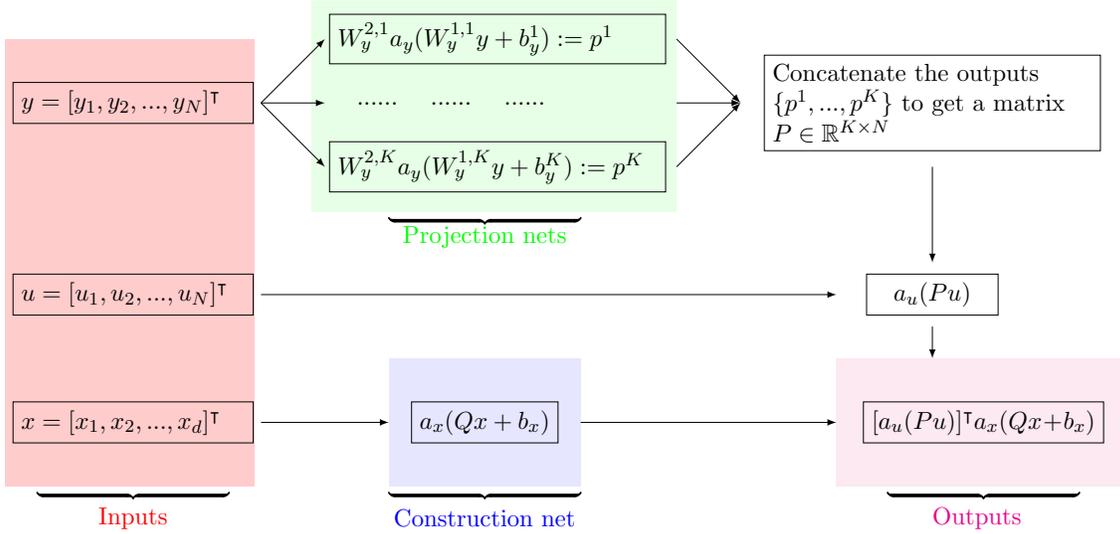

\textbf{Full Formulation}: BelNet trains both functions $p$ and $q$ directly. Intuitively, we project $u\in V$ into space spanned by trainable basis $p$, thus the name \textit{basis enhanced learning}. The network also constructs the trainable basis $q$ for the $x$-dependence, which we denote as \textit{construction net}. The entire architecture is displayed in Figure \ref{fig_network_s1}.

\subsection{Connection to FNO}
\label{sec_intuition}
The FNO form is motivated by iterating a block of linear layers in the Fourier domain, which is essentially a convolution representation for linear operators. To explain the connection between BelNet and FNO, let's consider the Fourier integral operator (Equation 4 in \cite{li2020fourier}). Essentially, the basis elements are trained in the full FNO form. In particular, consider the discrete Fourier transformation (DFT) of two functions $u$ and $h$, denoted as:
\begin{align}
    &\mathcal{F}(u)(k) = \sum_n u(n)b(n, k), \nonumber\\
    &\mathcal{F}(h)(k) = \sum_n h(n)b(n, k),
    \label{DFT}
\end{align}
where $b(n, k) = e^{-i2\pi nk}$ is the Fourier basis where $k$ denotes the wave number. By the discrete convolution theorem, we have:
\begin{align}
    u*h(m) = \sum_n u(n)h(m-n) = \mathcal{F}^{-1}\big(\mathcal{F}(u*h)\big) = \mathcal{F}^{-1}\left(\mathcal{F}(u)(k)\,\mathcal{F}(h)(k)\right).
\label{convolution}
\end{align}
Substituting Equation~\eqref{DFT} into Equation~\eqref{convolution}, it follows that,
\begin{align*}
    u*h(m) & = \sum_k\left( \sum_n u(n)b(n, k)\right) \left( \sum_n h(n)b(n, k) \right)w(k, m) = \sum_k w(k, m) \sum_n r(n, k) u(n) 
\end{align*}
where $w(k, m) = e^{i2\pi km}$ is the inverse DFT basis, and:
\begin{align*}
    r(n, k) = b(n, k)\sum_lh(l)b(l,k).
\end{align*}
Thus the linear FNO trains the following neural operator:
\begin{align*}
    G(u)(x) \approx \mathcal{FNO}(u)(m) = \sum_n u(n)h(m-n) = \sum_k w(k, m) \sum_n r(n, k) u(n)
\end{align*}
where $h$ is the kernel function. As a result, the linear Fourier integral operator is a special case of Equation~\eqref{form_linear}, where $b$ and $w$ play the role of $q$ and $p$ in Equation~\eqref{form_linear}. Note that the full FNO form includes an additional skip-connection and activation layer, which makes FNO a nonlinear approximation. Note that even in the full form, FNO uses the FFT and thus relies on the Cartesian domain with a lattice
grid mesh \cite{lu2022comprehensive}. BelNet does not have this restriction.

\section{Numerical Examples}
\label{sec_numerical}
We focus on four test PDEs: Burgers' equation, often used to compare neural operators, and three other challenging multiscale problems. While of great interest in the applied community \cite{hou1997multiscale, efendiev2013generalized}, multiscale and high-contrast problems pose several difficulties for neural operators. More precisely, one often needs a fine enough resolution to capture the  small-scale features, which leads to intractable training problems. In addition, observed data is often mesoscale or macroscale, and thus one needs to infer the homogenized solution operator. Our numerical experiments focus on predicting the solution based on data scarce observations. 

For each example, we test BelNet, stacked DON, and unstacked DON. In each case, we present the best results over the tests. We use two variants of BelNet; the first adopts the discretization used in DON experiments, which is referred to as fix-BelNet (since the mesh is the same for all samples). 
The second variant uses a randomized discretization for each sample; that is, each input function has a randomized mesh. This is referred to as free-BelNet. As a benchmark, we also refer to established methods \cite{efendiev2013generalized,hou1997multiscale,chetverushkin2021computational,chung2022computational,chung2018constraint} to solve multiscale parametric PDEs.

\subsection{Viscous Burgers' Equation}
\label{sec_vburgers}
Consider the viscous Burgers' equation with periodic boundary conditions: 
\begin{align*}
    &\frac{\partial u_s}{\partial t} + \frac{1}{2}\frac{\partial (u^2_s)}{\partial x} = \alpha \frac{\partial^2 u_s}{\partial x^2},\hspace{0.2em} x\in[0, 2\pi], \hspace{0.2em} t\in[0, 0.3]\\
    &u_s(x, 0) = u^0_s(x),\\
    &u_s(0, t) = u_s(2\pi, t),
\end{align*}
where $u^0_s(x)$ is the initial condition that depends on the parameter $s$ and the viscosity is set to $\alpha = 0.1$.
We will consider the operator that maps from the initial condition to the terminal solution at $t = 0.3$.

\begin{figure}[H]
\centering
\includegraphics[scale = 0.55]{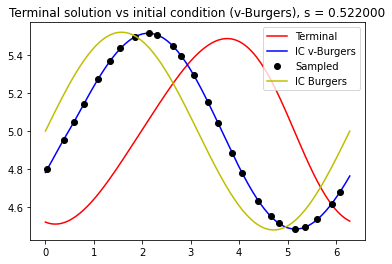}
\includegraphics[scale = 0.55]{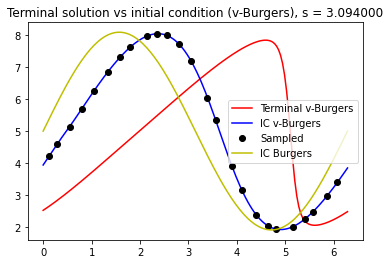}
\caption{Plots of two solutions to the viscous Burgers' equation with our initialization procedure. Note that each example's sampling points (black dots) for the initial condition differ. The yellow curves are used to generate the initial conditions for the model problem (viscous Burgers' equation). The initial conditions for the viscous Burgers' equation are displayed in blue.}
\label{pic_vburgers_info}
\end{figure}

\textbf{Training Data:} To generate the initial conditions, we first compute the short-time solution ($t=0.1$) to Burgers' equation. We use the periodic boundary conditions, set the viscosity to zero, and use the initial condition $s\sin(x)$ where $s\in [0, 4]$. The solution of the system at $t = 0.1$ is then used as the initial condition $u^0_s$; see the yellow and blue curves in Figure \ref{pic_vburgers_info}. This is used to generate more variability between initial samples for the training phase since different values of $s$ will lead to different levels of sharpness in the slopes for the initial data, even in the short time interval used.  
The mesh for the input data is as follows. Each initial condition (input function) has $25$ sensors, and we used a total of $200$ initial conditions for training. For each initial condition, the true system is evolved up to time $t=0.3$, and a total of $5$ time-stamps are collected (the terminal time is not included). Therefore the space-time mesh contains $25$-by-$5$ total sample locations for each initial condition.

\textbf{Testing:} To test the neural operators' ability to extrapolate future states,  we do not include the solution at the terminal time $t=0.3$ in the training dataset. For testing, we use solutions from $500$ initial conditions and test each neural operator on the solution at the terminal time with a finer mesh of $151$ grid points. 
We present the relative errors in Table~\ref{vburgers_results}. We compare BelNet with the (vanilla) DON, where BelNet uses either fixed sensor points matching the discretization in DON or ``free'' sensors that are independent of each sample. With a similar amount of trainable parameters listed in Table~\ref{vburgers_results}, BelNet performs similarly or better in terms of the relative prediction error over the 500 testing sequences.

\begin{table}[H]
\centering
\begin{tabular}{||c c c c c||} 
\hline
 & DON & fix-BelNet  & free-BelNet & FNO\\ [0.5ex] 
\hline
Relative error (\%) & $1.17$ & $0.59$ &  $1.07$ & $1.61$ \\ [0.5ex]
\hline
Total parameters (in K) & $103.7$ & $102$ &  $102$ & $287.425$\\ [0.5ex]
\hline
\end{tabular}
\caption{Viscous Burger's equation example mean relative error of $500$ testing solutions of four experiments. The FNO structure is from \cite{li2020fourier}.}
\label{vburgers_results}
\end{table}

\textbf{BelNet Hyperparameters}: The number of sensors (for the initial condition) is set to $N = 25$. For the projection net, ``$u$'' is connected with $K = 10$ independent subnetworks with identical structures. 
Specifically, each sub-network is fully connected and has the structure $[25, 100, 25] \quad (i.e.\ N = 25, N_1 = 100)$. Here $25$, $100$, and $25$ are the input dimension, the width of the first hidden layer (the number of neurons in the first layer), and the output dimension, respectively.
This network representation routine will be used in the rest of the manuscript. In addition, bias and activation (Tanh) are added to the first layer according to the formulation (\ref{intro_formulation3}).
\textcolor{black}{To achieve better accuracy, we replace $a_u$ in Figure \ref{fig_network_s1} with a ReLU-activated trainable layer of size $10\times 10$.} More precisely, we use a ReLU-activated trainable matrix of size $10\times 10$ to replace the activation function $a_u$ in Figure \ref{fig_network_s1}. 
For the construction net, the input dimension equals $2$ ($1$ in space and $1$ in time). Hence it is connected with a fully connected network, which has the structure $[2, 100, 100, 100, 10] \quad (i.e.\ K = 10)$. Each layer except the last layer has a bias vector and is activated by the Tanh function.

\subsection{Multiscale Elliptic Equation}
\begin{figure}[t]
\centering
\includegraphics[scale = 0.5]{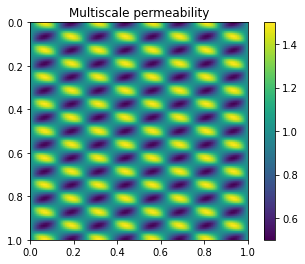}
\includegraphics[scale = 0.4]{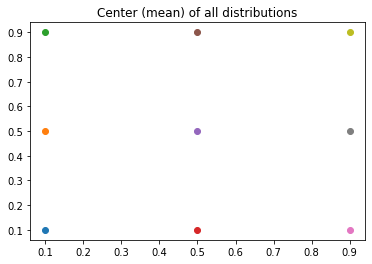}\\
\includegraphics[scale = 0.5]{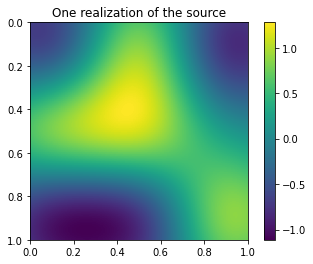}\hspace{0.1em}
\includegraphics[scale = 0.5]{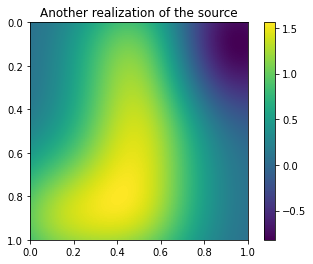}
\caption{Plots (from top left to bottom right): the multiscale permeability, centers of all normal distributions, one realization (sample) of the source, a second realization of the source. }
\label{pic_elliptic_info}
\end{figure}

We consider the following multiscale elliptic equation:
\begin{align*}
    & -\nabla \cdot (\kappa \nabla u) = f(x, s),\ \  x\in s = [0, 1]^2, \ \ s\in D\\
    & u(x) = 0,\ \  x\in\partial \Omega,
\end{align*}
where $s$ is the parameter. The multiscale permeability is denoted by $\kappa$ and an example is shown in Figure \ref{pic_elliptic_info}, specifically:
\begin{align*}
    \kappa(x_1, x_2) = 1 + \frac{\sin(2\pi\frac{x_1}{\epsilon_1})\cos(2\pi\frac{x_2}{\epsilon_2}) }{2+\cos(2\pi\frac{x_1}{\epsilon_1})\sin(2\pi\frac{x_2}{\epsilon_2})}
    +  \frac{\sin(2\pi\frac{x_1}{\epsilon_1})\cos(2\pi\frac{x_2}{\epsilon_3}) }{2+\cos(2\pi\frac{x_1}{\epsilon_1})\sin(2\pi\frac{x_2}{\epsilon_3})},
\end{align*}
where the parameters are set to $\epsilon_1 = 1/4$, $\epsilon_2 = 1/8$ and $\epsilon_3 = 1/16$.
The multiscale nature of this problem comes from the permeability $\kappa(x_1, x_2)$, which has multiple frequency scales. General machine learning-based approaches have difficulty approximating such spatial systems due to the multi-frequency behavior of the solution \cite{leung2022nh}.
The source function, denoted by $f$, is parameterized by taking the weighted sum of $9$ Gaussian distributions whose centers are uniformly distributed in the domain. A demonstration is presented in Figure \ref{pic_elliptic_info}. We sample the weights for the Gaussians in the source function such that they are uniformly from $[-1, 1]$. The training problem is to learn the map from the source $f$ to the solution of the elliptic equation, $u$. 

\textbf{Training/Testing Data:} Using the parametric form for $f$, we generate $80$ different source functions for training. For each solution $u$ (corresponding to one source $f$), we choose $19\times 19$ grid points in $\Omega$ uniformly, leading to a total of $80\times 19 \times 19$ training ``samples''. For testing, we generate $100$ testing permeability coefficients and compute the predicted solution on a fine mesh using $100\times 100$ grid points. \textcolor{black}{Each parameter $s$ will determine an input function $f(x, s)$. To measure the input function, we uniform (randomly)  placed $10\times 10$ sensors in the input function domain $\Omega$, i.e. $N = 100$. More precisely, we create a uniform coarse mesh of size $10\times 10$ and draw one sensor uniformly from each coarse element. }

\textbf{Results and Comparison:} We compare the results with provably convergent numerical methods in order to better evaluate our approach. This system is a classical multiscale problem, a standard finite element solver would be too slow since one would need a fine discretization level (to capture the multiscale features and obtain sufficient accuracy). Therefore, we compare our approach to the multiscale finite element method (MsFEM) \cite{hou1997multiscale, efendiev2013generalized}, which is an efficient and accurate method for this problem. The fine-scale mesh size is set to
$100\times 100$, and the coarse element size is set to $5\times 5$. Thus, we have $19\times 19$ degrees of freedom, which is the same as the number of solution points per training source used in the operator learning approaches in this example. 
The results are summarized in Table \ref{elliptic_results}.

\begin{table}[H]
\centering
\begin{tabular}{||c c c c c||} 
\hline
& DON (stack) & fix-BelNet  & free-BelNet & MsFEM \\ [0.5ex] 
\hline
Euclidean error (\%) & $16.44$ & $14.43$ & $9.66$ & $2.28$\\ [0.5ex]
\hline
Parameters count (in K) & $319.4$ & $312.4$ &  $312.4$ & NA\\ [0.5ex]
\hline
\end{tabular}
\caption{Elliptic example mean relative error of $100$ testing solutions of three experiments. 
MsFEM adopts $19\times 19$ DOFs while the learning methods are trained on the same $19\times 19$ mesh. 
}
\label{elliptic_results}
\end{table}

\textbf{BelNet Hyperparameters:}
As before, we set $N = 100$. For the projection net, $u$ is connected with $K = 10$ independent networks with the structure $[200, 100, 100] \quad (i.e.\ N = 100, N_1 = 100)$. In addition, the bias is added to the first layer, and the first layer is activated by the ReLU function. Since the problem is linear in $f$, to accelerate the training, we removed the nonlinear activation $a_u$. For the construction net, the input dimension is equal to $2$ and is a fully connected network of a structure $[2, 100, 10] \quad (i.e.\ K = 10)$. Each layer except the last layer has a bias vector and is activated by the ReLU function.

\subsection{Radiative Transfer Equation (RTE) with High Contrast Scatter}
\begin{figure}[H]
\centering
\includegraphics[scale = 0.35]{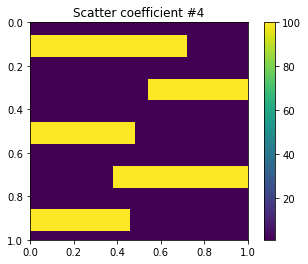}
\includegraphics[scale = 0.35]{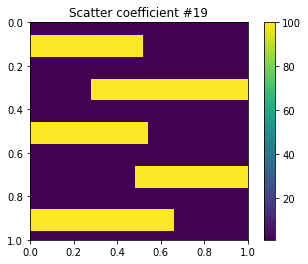}
\includegraphics[scale = 0.35]{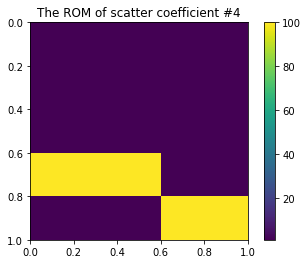}
\includegraphics[scale = 0.35]{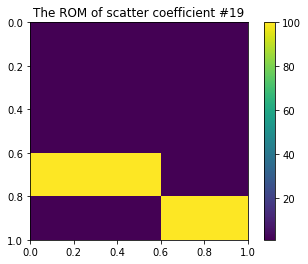}
\caption{Left two images: two scatter coefficients. Right two images: the corresponding ROMs of two scatter coefficients. The ROM is obtained by sampling $\sigma(x, \omega)$ in a uniform mesh $5\times 5$ in the domain, and two scatter coefficients are sampled by the same mesh. }
\label{exp_demon_sampling}
\end{figure}

We consider the radiative transfer equation (RTE) \cite{newton2020diffusive, li2019diffusion, chung2020generalized} with a high contrast scatter coefficient $\sigma(x, \omega)$ (see Figure \ref{exp_demon_sampling}).
Specifically, the `high contrast' aspect comes from strong scattering in the channels. We assume that the channels are moving and the channels are characterized implicitly by the parameter $\omega$. The goal is to learn the mapping from $\sigma(x, \omega)$ with moving channels to the solution, $I(x,s)$, of the RTE problem:
\begin{align*}
 s \cdot \nabla  I(x,s) = \frac{\sigma(x, \omega)}{\epsilon} \bigg( \int_{\mathcal{S}^{n-1}} I(x,s') ds' - I(x, s)\bigg) \quad \forall x \in D, s \in \mathcal{S}^{n-1}.
\end{align*}
Here $s$ is a vector on the unit sphere, and $n$ is the dimension of the problem. In our experiments, we considered $n=2$ and thus $\mathcal{S}^{n-1}=\mathcal{S}^{1}$ is the unit circle. 
In addition, we set $\epsilon = 1$ and $D = [0, 1]^2$.
To provide closure to the model, we introduce the Dirichlet boundary conditions $I(x,s) = I_{\text{in}}$ used for entrant directions $s \cdot \textbf{n} <0$, i.e. on $\Gamma ^{-} := \{ (x,s) \in \partial D \times \mathcal{S}^{n-1}: s \cdot \textbf{n} <0 \}$. Here, $\textbf{n}$ is the unit outward normal vector field at $x \in \partial D$. The condition is written as 
\begin{align*}
I = I_{\text{in}} (x,s) \quad \text{ for all } (x,s) \in \Gamma^-.
\end{align*}
In our examples, the top, bottom, and right boundaries have zero incoming boundary conditions, while we assume the left boundary has non-zero flow injected into the domain.

We chose the multiscale RTE with high contrast channels since it is complicated to solve numerically \cite{chung2020generalized} and a challenging problem for learning-based approaches. In practice, it is impossible to measure $\sigma(x, \omega)$ at all $x\in D$. Instead, one samples $\sigma(x)$ on a mesh $D'$ computes a reduced order model $\hat{\sigma}(x, \omega)$. However, due to the nature of the scattering coefficient, i.e. sparse localized regions, the ROM for two different scatter coefficients may appear to be similar even when the original structures differ. This is illustrated in Figure \ref{exp_demon_sampling} and is often the case when one uses the same grids to construct the ROM. Consequently, a mesh-free version of operator learning should be used. 

\textbf{Training/Testing Data:} We generate $60$ scattering coefficients and compute the corresponding solutions $I(x,s)$. For each solution, we use $25\times 25$ points in $D$ uniformly sampled with $8$ velocities, which gives a total of $60\times 625 \times 8$ training samples. For testing, we used $100$ scatter coefficients and computed the corresponding solutions on a finer mesh of $50\times50\times 8$, where the last dimension is the velocity. \textcolor{black}{Each $\omega$ determines one scatter coefficient $\sigma(x, \omega)$. To measure each input function $\sigma(x, \omega)$, we draw $25$ sensors randomly uniformly from $D$, so $N = 25$. More precisely, we create a uniform discretization of size $5\times 5$ and draw one sensor randomly from each element in the mesh.}

\begin{table}[H]
\centering
\begin{tabular}{||c | c| c| c | c||} 
\hline
\backslashbox{Model}{CR} & $2$ & $10$  & $50$ & $100$ \\ [0.5ex] 
\hline
fix-BelNet (279K) & $3.68 $  & $8.36 $ & $13.70$ &  $16.40$\\ [0.5ex]
\hline
free-BelNet (279K)& $3.56$ & $7.28$ & $12.30$ & $13.31$\\ [0.5ex]
\hline
DON  (294K)& $4.39$ & $15.12$ & $18.30$ & $19.18$\\ [0.5ex]
\hline
\end{tabular}
\caption{RTE mean relative error of $100$ testing solutions for different contrast ratios (CR). Fix-BelNet and free-BelNet share the structure and have 279K trainable parameters, DON has $294K$ parameters. }
\label{rte_cr_table}
\end{table}

\textbf{Results and Comparison:} Table \ref{rte_cr_table} summarizes the results. As expected, in all three models the performance deteriorates as the contrast ratio (CR) increases. However, it is worth noting that the grid-free version of BelNet seems to produce a small relative error with a smaller jump between various CR values.
This may be due to underlying sampling of the grid points, namely since the grids are randomly generated free-BelNet sees more variations between inputs.
Since fix-BelNet has a lower relative error than the DON, this indicates that incorporating mesh coordinate information indeed improves the model's accuracy.

\textbf{BelNet Hyperparameters:} We set $N = 25$, the projection net's parameters are $K = 30$ with fully connected structure of size $[50, 100, 25] \quad (i.e.\ N = 25, N_1 = 100)$. As before, the bias is added to the first layer and the first layer is activated by the ReLU function. For the construction net, the input dimension is equal to $4$ ($2$ in space and $2$ in velocity), hence it has the structure $[4, 100, 100, 30] \quad (i.e. \ K = 30)$. Each layer except the last layer of the construction net has a bias vector and is activated by the ReLU function.

\subsection{Richard's Equation}
\begin{figure}[H]
\centering
\includegraphics[scale = 0.5]{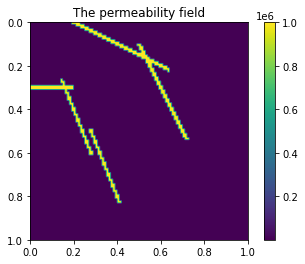}
\includegraphics[scale = 0.5]{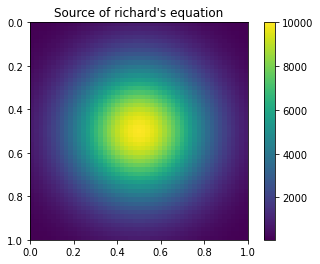}
\caption{Left: the permeability field $\kappa(x)$ used in Richard's equation. The permeability in the yellow channels is equal to $1$ while the permeability in the background is equal to $10^6$. Right: the source.}
\label{kappa_richard}
\end{figure}

The last example is the nonlinear Richard's equation \cite{efendiev2022efficient, efendiev2022hybrid}:
\begin{align*}
    & u_t(x, t) - \nabla \cdot e^{\kappa(x) u(x, t)} \nabla u(x, t) = f(x), \ \ x\in\Omega = [0, 1]^2, \ \ t\in[0, T],\\
    & \frac{\partial u(x, t)}{\partial n} = 0, \ \ x\in\partial \Omega,
\end{align*}
with a high-contrast multiscale contrast permeability  $\kappa (x)$\cite{chung2022contrast, efendiev2022efficient, efendiev2022hybrid} (see Figure \ref{kappa_richard}). The source is defined as $$f(x, y) = 10^4 \, \exp\left(\frac{-(x-0.5)^2-(y-0.5)^2}{0.1}\right).$$

\textbf{Training/Testing Data:} We map the initial condition to the terminal solution at time $t = 2\times 10^{-5}$.
The initial condition is defined as the weighted sum of $9$ Gaussian distributions, where the weights are uniformly random from $[-2, 2]$. The training set is generated from $60$ initial conditions, where for each initial condition, we observe the solution on a $40\times 40\times 7$ mesh, i.e. $40\times 40$ spatial points and $7$ time-snapshots. For testing, we test the model on $100$ solutions only at the terminal time. For the test solutions, we use a (finer) spatial mesh of size $50\times 50$. To measure the initial condition, free-BelNet randomly draws $25$ points for each sample, while the DON chooses $25$ fixed points for all training and testing samples; in addition, we also test fix-BelNet with the same fixed points used in DON.

\textbf{Results and Comparison:} We compare our results to DON and the solution generated by the generalized multiscale finite element method (GMsFEM) \cite{chetverushkin2021computational, chung2022computational}. For GMsFEM, the fine mesh size is set to $50\times 50$, the coarse element size is $2\times 2$, and we use the first $3$ GMsFEM basis. As a result, GMsFEM has $1728$ degrees of freedom, which is close to the $40\times 40$ training mesh points used in the operator learning models. The results are summarized in Table \ref{40nl_results}.

\begin{table}[H]
\centering
\begin{tabular}{||c c c c c c||} 
\hline
 & DON (unstack) & DON (stack) & fix-BelNet  & free-BelNet & GMsFEM \\ [0.5ex] 
\hline
Euclidean error (\%) & $11.98$ & $21.18$ & $7.87$ &  $7.64$ & $32.43$\\ [0.5ex]
\hline
Parameters count (in K) & $103$ & $134$ & $127$ &  $127$ & NA \\ [0.5ex]
\hline
\end{tabular}

\caption{Richard's equation mean relative error of $100$ testing solutions of three experiments. GMsFEM adopts $1,728$ DOFs while the learning methods are trained on $40\times 40$ mesh.}
\label{40nl_results}
\end{table}

\textbf{BelNet Hyperparameters:} We set $N = 25$. For the projection net, we use $K = 10$ and the fully connected network has the structure $[50, 100, 25] \quad (i.e.\  N = 25, N_1 = 100)$, with bias in the first layer and with the Tanh activation function. 
In addition, we replace $a_u$ in Figure \ref{fig_network_s1} with a Tanh-activated trainable layer of size $10\times 10$, which improves accuracy. For the construction net, we use a $[3, 100, 100, 10] \quad (i.e.\  K = 10)$ fully connected neural network with bias and Tahn activation on each layer except the last layer.

\section{Conclusions and Future Works}
We proposed  a mesh-free operator learning technique called BelNet. In particular, a method that does not depend on the grid locations for the input function and allows for the domain and discretization of the inputs and outputs to differ. Mesh-free approaches make it possible to solve challenging PDE problems, such as multiscale systems where functions can be aliased more often for fixed grids. 
While the form of BelNet is motivated by the convolution representation for linear operators related to PDE,  we provide a nonlinear generalization using additional activations. Through several detailed examples, we show that BelNet can outperform standard and recent approaches on important parametric PDEs. 

There are several future directions worth investigating. First, the universal approximation theorem for BelNet should follow from a generalization of \cite{chen1995universal}, but one needs to show that the approximation does not require fixed or similar grids. Secondly, in order to improve the computational and storage cost, we plan to optimize the parameter count in the projection net. In particular, leveraging the success of random feature projections, we have seen some preliminary improvements to the efficiency of BelNet. Lastly, while we focused on the applications to multiscale PDE,  BelNet can also be applied to complicated input controls for dynamical systems with controls.

\section{Acknowledgement}
Z. Zhang was supported in part by AFOSR MURI FA9550-21-1-0084. H. Schaeffer was supported in part by AFOSR MURI FA9550-21-1-0084.

\bibliographystyle{abbrv}
\bibliography{references}
\end{document}